\documentclass[11pt]{amsart}
\usepackage{graphicx,amssymb,amsmath,amsthm}
\usepackage{enumerate}
\usepackage{comment,cite,color}
\usepackage{cite,color}
\usepackage{algorithmic}
\usepackage{mathrsfs}
\usepackage[ruled]{algorithm}
\usepackage{epsfig}
\usepackage{lscape}

\usepackage{algorithm}  
\usepackage{algorithmic} 

\newtheorem{theorem}{Theorem}[section]

\newcommand{\PP}{\mathbb{P}}

\newcommand{\innerp}[1]{\langle {#1} \rangle}

\newcommand{\abs}[1]{\lvert#1\rvert}

\newtheorem{remark}{Remark}

\renewcommand{\SS}{{\mathcal S}}

\newcommand{\R}{{\mathbb R}}
\newcommand{\Z}{{\mathbb Z}}

\newcommand{\C}{{\mathbb C}}
\newcommand{\CG}{{\mathcal G}}
\newcommand{\Q}{{\mathbb Q}}

\renewcommand{\eqref}[1]{(\ref{#1})}

\newcommand{\ML}{{\mathcal L}}
\newcommand{\MLP}{{\mathcal{PL}}}

\renewcommand{\H}{{\mathbb H}}

\newtheorem{prop}{Proposition}[section]
\newtheorem{lem}[prop]{Lemma}

\newtheorem{coro}[prop]{Corollary}

\newtheorem{conjecture}[prop]{Conjecture}

\usepackage[top=0.6in,bottom=0.6in,left=0.7in,right=0.7in]{geometry}
\date{}

\begin{document}
\bibliographystyle{plain}
\title{The minimal measurement number for  low-rank matrices recovery}

\author{Zhiqiang Xu}
\thanks{  Zhiqiang Xu was supported  by NSFC grant (11171336, 11422113, 11021101, 11331012) and by National Basic Research Program of China (973 Program 2015CB856000).}
\address{LSEC, Inst.~Comp.~Math., Academy of
Mathematics and System Science,  Chinese Academy of Sciences, Beijing, 100091, China}
\email{xuzq@lsec.cc.ac.cn}
\begin{abstract}
The paper presents several results that address a fundamental question in low-rank matrices recovery: how many measurements
are needed to recover low rank matrices? We begin by investigating the complex matrices case and show that $4nr-4r^2$ generic measurements are both necessary and sufficient for the recovery of  rank-$r$ matrices in $\C^{n\times n}$ by algebraic tools developed in \cite{CEHV}. Thus, we confirm a conjecture which is raised by Eldar, Needell and Plan for the complex case. We next consider the real case and prove that the bound  $4nr-4r^2$ is tight provided $n=2^k+r, k\in \Z_+$.
Motivated by Vinzant's work,  we construct $11$ matrices in $\R^{4\times 4}$ by computer random search and prove they define injective measurements on rank-$1$ matrices in $\R^{4\times 4}$. This disproves the conjecture raised by Eldar, Needell and Plan for the real case. Finally, we use the results in
this paper to investigate the phase retrieval by projection and  show fewer than $2n-1$ orthogonal projections are possible for the recovery
 of $x\in \R^n$ from the norm of them, which gives a negative answer for a question raised in \cite{phaseproj1}.
\end{abstract}

\maketitle

\bigskip \medskip

\section{Introduction}
\subsection{Problem setup}
The problem of low-rank matrix recovery attracted many attention recently  since it is widely used in
 image processing,  system identification and control, Euclidean embedding, and recommender systems.
 Suppose that the matrix $Q\in \H^{n\times n}$ with ${\rm rank}(Q)\leq r$, where $\H$ is either $\R$ or $\C$. The information we gather about $Q$ is
\[
b_j:=\innerp{A_j,Q}:={\rm trace}(A_j^*Q), \quad j=1,\ldots,m
\]
where $A_j\in \H^{n\times n}, j=1,\ldots, m$. The aim of the low-rank matrix recovery is to recover $Q$ from ${\bf b}=[b_1,\ldots,b_m]\in \H^m$.

For a given ${\mathcal A}:=\{A_1,\ldots,A_m\}\subset \H^{n\times n}$, we define the map ${\bf M}_{\mathcal A}:\H^{n\times n}\rightarrow \H^m$
by
\[
{\bf M}_{\mathcal A}(Q)=[b_1,\ldots,b_m].
\]
Set
\[
\ML_r^\H:=\{X\in \H^{n\times n}: {\rm rank}(X)\leq r\}.
\]
 We say the matrices set  ${\mathcal A}:=\{A_1,\ldots,A_m\}$ has the low-rank matrix recovery property for $\ML_r^\H$ if the map ${\bf M}_{\mathcal A}$ is injective on $\ML_r^\H$.
Naturally, we are interested in the minimal $m$ for which the map ${\bf M}_{\mathcal A}$ is injective on $\ML_r^\H$.

There are many convex programs for the recovery of the low-rank matrix $Q$ from  ${\bf M}_{\mathcal A}(Q)$. A well-known one is nuclear-norm minimization which requires $m=Cnr$ random linear measurements for the recovery of rank-$r$ matrices in $\H^{n\times n}$ \cite{LMCR,LMPlan,LMRe,LMshen}.
 Despite many literatures on  low-rank matrices recovery, there remains a fundamental lack of understanding
about the theoretical limit of the number of the cardinality of ${\mathcal A}$ which
has the low-rank matrix recovery property for $\ML_r^\H$. This paper focusses on the problem of the minimal measurements number for
 the recovery of low-rank matrix. We state the problem as follows:
\begin{enumerate}[{\bf Problem }1]
 \item What is the minimal measurement number $m$ for which there exists ${\mathcal A}=\{A_1,\ldots,A_m\}\subset \H^{n\times n}$ so that ${\bf M}_{\mathcal A}$ is injective on $\ML_r^\H$?
\end{enumerate}
The aim of this paper is to  addresses Problem  1 under many different settings.

\subsection{Related work}
A related problem to low-rank matrices recovery is {\em phase retrieval}, which is to recover  the rank-one matrix $xx^*\in \H^{n\times n}$ from
the measurements $\abs{\innerp{\phi_j,x}}^2=\innerp{\phi_j\phi_j^*,xx^* }, j=1,\ldots,m,$ where $\phi_j\in \H^n$ and $x\in \H^n$. In the context of phase retrieval,
one is interested in the minimal measurement number $m$ for which the map ${\bf M}_{\Phi}$ is injective on $\SS_1^\H$ where $\Phi:=\{\phi_1\phi_1^*,\ldots,\phi_m\phi_m^*\}$ and $\SS_r^\H:=\{X\in \H^{n\times n}: {\rm rank}(X)\leq r, X^*=X\}, r\in \Z$.
It is known that in the real case $\H=\R$ one needs at least $m\geq 2n-1$ vectors so that ${\bf M}_{\Phi}$ is injective on $S^\R_1$ \cite{BCE}. For the
complex case $\H=\C$, the same problem remain open.
 Balan, Casazza and Edidin  in \cite{BCE}  show that
${\bf M}_{\mathcal A}$  is injective on $S^\C_1$ if $m\geq 4n-2$ and $\phi_1,\ldots,\phi_m$ are generic vectors in $\C^n$.
In \cite{BCMN}, Bandeira, Cahill, Mixon and Nelson conjectured the following (a) if $m<4n-4$ then ${\bf M}_\Phi$ is not injective on $\SS_1^\C$;
(b) if $m\geq 4n-4$ then ${\bf M}_\Phi$ is  injective on $\SS_1^\C$ for generic vectors $\phi_j, j=1,\ldots,m$.
The part (b) of the conjecture is proved by Conca,  Edidin,   Hering and  Vinzant in \cite{CEHV}    by employing  algebraic tools.
They also confirm  part (a) for the case where $n$ is in the form of $2^k+1, k\in \Z$.  Recently, in \cite{small}, a counterexample is presented disproving part (a) of this conjecture. In fact, \cite{small} gives $11=4n-5<4n-4$ vectors ${\phi_1},\ldots,{\phi_{11}}\in \C^{4}$ and prove that ${\bf M}_{ \Phi}$ is injective on $\SS_1^\C$ by algebraic computation where ${ \Phi}=\{{\phi}_1{\phi}_1^*,\ldots, {\phi}_{11}{\phi}_{11}^* \}$.

In context of low-rank matrix recovery, it  is Eldar, Needell and Plan \cite{uniq} that show that
$m\geq 4nr-4r^2$ Gaussian matrices $A_1,\ldots,A_m$ has low-rank matrix recovery property  for $\ML_r^\H$ with probability 1 (see also \cite{His1,His2,His3}) provided $r\leq n/2$. Naturally, one may be interested in whether the number $4nr-4r^2$ is tight. In \cite{uniq}, the authors  made the following conjecture:

\begin{conjecture}\cite{uniq}
{ If $m<4nr-4r^2$ then ${\bf M}_{\mathcal A}$ is not injective on $\ML_r^\H$.}
\end{conjecture}

\subsection{Our contribution}

The aim of this paper is to address Problem 1 by employing algebraic tools which are developed in \cite{CEHV,small}. In Section 2, we consider the
case $\H=\C$ and prove that ${\bf M}_{\mathcal A}$ is injective on $\ML^\C_r$ if ${\mathcal A}$ contains $m\geq 4nr-4r^2$ generic matrices $[A_1,\ldots,A_m]\in \C^{mn^2}$. Compared to the results of \cite{uniq}, we do not require the $m$ matrices are Gaussian matrices. Hence,
our result does not suffer from probabilistic qualifiers on the injective  (E.g., injective ``with probability 1").
We also show that the bound $4nr-4r^2$ is tight which means ${\bf M}_{\mathcal A}$ is not injective
on $\ML^\C_r$ provided $m<4nr-4r^2$ and hence confirm Conjecture 1.1 for the complex case, i.e., $\H=\C$.   We turn to the real case in Section 3 and prove the bound  $4nr-4r^2$ is tight provided $n$ is in the form of $2^k+r$.
Inspired by the work of \cite{small}, we use computer random search to construct a counterexample for the case $n=4, r=1$. In fact, we present $11=4n-5$ matrices $A_1,\ldots,A_{11}\in \R^{4\times 4}$ and prove ${\bf M}_{\mathcal A}$ is injective on $\ML_1^\R$ using Vinzant's test with disproving  Conjecture 1.1 for the real case. We next consider the recovery of the symmetric matrix and  investigate the minimal measurement number $m$ for which there exists ${\mathcal A}=\{A_1,\ldots,A_m\}\subset \R^{n\times n}$ so that ${\bf M}_{\mathcal A}$ is injective on $\SS_r^\R$. In Section 4, we apply the results to study the phase retrieval by projection.
Set $W_j:={\rm span}\{u_{j,1},\ldots,u_{j,d_j}\}\subset \R^n$ and $P_j:\R^n\rightarrow W_j$ is an orthogonal projection.
Following \cite{phaseproj,phaseproj1},  phase retrieval by projection is to recover $x\in \R^n$  up to a unimodular constant  from $\|P_jx\|_2$ .
We say that $\{W_j\}_{j=1}^m$ yields phase retrieval if for all $x,y\in \R^n$ satisfying $\|P_jx\|=\|P_jy\|$ for all $j=1,\ldots,m$ then $x=\pm y$.
In \cite{phaseproj}, Cahill, Casazza, Peterson and  Woodland proved that $2n-1$ projections are enough for phase retrieval. Particularly, they showed that phase retrieval can be done in $\R^n$ with $2n-1$ subspaces each of any dimension less than $n-1$.
A question  is also raised in \cite{phaseproj1} which states {\em can phase
retrieval be done in $\R^n$ with fewer than $2n-1$ projections?} Using the results in this paper, we present a positive answer for the question provided  $n$ is in the form of $2^k+1$. We also give a negative answer for the case $n=4$ by constructing $6=2n-2$ subspaces
$W_1,\ldots,W_6\subset \R^4$ and prove they are phase retrieval by computational algebra.

\section{The recovery of complex low rank matrices }

We first recall the following lemma
\begin{lem}\label{le:uniq}\cite{uniq}
Suppose that $r\leq n/2$. The map ${\bf M}_{\mathcal A}$ is not injective on $\ML^\H_r$ if and only if there is a nonzero $Q\in \ML^\H_{2r}$ for which
\[
{\bf M}_{\mathcal A}(Q)=0.
\]
\end{lem}

According to Lemma \ref{le:uniq}, the set $\ML^\C_{r}$ plays an important role in the investigation of the ${\bf M}_{\mathcal A}$.
Recall that ${\rm rank}(Q)\leq r$ is equivalent to the vanishing  of all $(r+1)\times (r+1)$ minors of $Q$. Hence, $\ML_r^\C$ is an affine variety in $\C^{n^2}$ with the dimension  $2nr-r^2$ \cite[Prop. 12.2]{alge}  and the degree
$d_{n,r}:=\prod_{i=0}^{n-r-1}\frac{(n+i)!\cdot i!}{(r+i)!\cdot (n-r+i)!}$ \cite[Ex. 19.10]{alge}. Note these $(r+1)\times (r+1)$ minors are  homogeneous polynomials in the entries of $Q$. Thus the projectivization of $\ML_r^\C$ is a projective variety in
$\PP(\C^{n^2})$ and it is called as determinant variety.  Throughout the paper, by  $m$ generic matrices in $\H^{n\times n}$  we mean  $[A_1,\ldots,A_m]$ corresponds to a point in a non-empty  Zariski open subset
of $\H^{mn^2}$ which is also open and dense in the Euclidean topology (see \cite[Section 2.2]{CEHV}). We next state the main
result of this section:

\begin{theorem}\label{th:comp}
Suppose that $r\leq n/2$.
Consider $m$ matrices ${\mathcal A}=\{A_1,\ldots,A_m\}\subset \C^{n\times n}$ and the mapping ${\bf M}_{\mathcal A}:\C^{n\times n}\rightarrow \C^m$. The following holds

(a) If $m\geq 4nr-4r^2$ then ${\bf M}_{\mathcal A}$ is injective on $\ML_r^\C$ for generic matrices $A_1,\ldots,A_m$.

(b) If $m<4nr-4r^2$, then ${\bf M}_{\mathcal A}$ is not injective on $\ML_r^\C$.

\end{theorem}
\begin{proof}
We use $\CG_{m,n}$ to denote the matrices set  $([A_1,\ldots,A_m],[Q])\in \PP(\underbrace{\C^{n\times n}\times \C^{n\times n}\times\cdots \times\C^{n\times n}}_m)\times \PP(\C^{n\times n})$
 which satisfies  the following property:
\[
{\rm rank}(Q)\leq 2r\quad \text{and } \quad \innerp{A_j,Q}=0, \quad \text{ for all } 1\leq j\leq m.
\]
Note that $\CG_{m,n}$ is defined by the vanish of homogeneous polynomials in the entries of $A_j$ and $Q$. Thus
 $\CG_{m,n}$ is a projective variety of $\PP(\underbrace{\C^{n\times n}\times \C^{n\times n}\times\cdots \times\C^{n\times n}}_m)\times \PP(\C^{n\times n})$.
We next consider the dimension of the projective complex variety $\CG_{m,n}$. We let $\pi_1$ and $\pi_2$ be projections onto the first
and the second coordinates, respectively, i.e.,
\[
\pi_1([A_1,\ldots,A_m],[Q])=[A_1,\ldots,A_m],\quad  \pi_2([A_1,\ldots,A_m],[Q])=Q.
\]
We claim that $\pi_2(\CG_{m,n})=\MLP_{2r}^\C$ where
\[
\MLP_{2r}^\C:=\{Q\in \PP(\C^{n\times n}): {\rm rank}(Q)\leq 2r\}.
\]
 Indeed, for any fixed $Q_0\in \MLP_{2r}^\C$, there exists a matrix $A_0\in \C^{n\times n }$ satisfying
$\innerp{A_0,Q_0}=0$ since $\innerp{A_0,Q_0}$ is a linear equation about the entries of $A_0$. This implies that $([A_0,\ldots,A_0],[Q_0]) \in \CG_{m,n}$ and $\pi_2([A_0,\ldots,A_0],[Q_0])=Q_0$. Thus we have $\pi_2(\CG_{m,n})=\MLP_{2r}^\C$.
Note that $\MLP_{2r}^\C\subset \PP(\C^{n\times n})$ is a projective variety with   ${\rm dim}(\pi_2(\CG_{m,n}))=4nr-4r^2-1$.

We next consider the dimension of the preimage $\pi_2^{-1}(Q_0)\subset \PP(\underbrace{\C^{n\times n}\times \cdots \times \C^{n\times n}}_m)$ for a fixed $Q_0\in \PP(\C^{n\times n})$. A simple observation is that
\begin{equation}\label{eq:plane}
\innerp{A_j,Q_0}=0
\end{equation}
defines a nonzero linear equation on the entries of $A_j$. For each $A_j$ the linear equation (\ref{eq:plane}) defines a hyperplane of dimension $n^2-1$ in $\C^{n^2} \cong \C^{n\times n}$. Hence, after projectivization, the preimage $\pi_2^{-1}(Q_0)$ has dimension $m(n^2-1)-1=mn^2-m-1$. Then,
according to \cite[Cor.11.13]{alge}
\[
{\rm dim}(\CG_{m,n})= {\rm dim}(\pi_2(\CG_{m,n}))+{\rm dim}(\pi_2^{-1}(Q_0))=(4nr-4r^2-1)+(mn^2-m-1)=mn^2+4nr-4r^2-m-2.
\]
If $m\geq 4nr-4r^2$, then
\begin{equation}\label{eq:dimpi1}
{\rm dim}(\pi_1(\CG_{m,n})) \leq {\rm dim}(\CG_{m,n})=mn^2+4nr-4r^2-m-2<mn^2-1.
\end{equation}
 Here, we use the result which states  the dimension of the projection  is less  or equal to the dimension of the original variety \cite[Cor.11.13]{alge}. Note the dimension of the
$\PP(\underbrace{\C^{n\times n}\times \cdots \times \C^{n\times n}}_m)$, which is the target of the projection $\pi_1$, is $mn^2-1$. The (\ref{eq:dimpi1}) shows the dimension of $\pi_1(\CG_{m,n})$ is strictly less than $mn^2-1$ provided $m\geq 4nr-4r^2$. This means the image of the projection $\pi_1$ lies in a hyper-surface which is defined by the vanish of some polynomials. We arrive at (a).

We next turn to (b).
For ${\mathcal A}=\{A_1,\ldots,A_m\}$,
we set
\[
Z_{\mathcal A} \,\,:=\,\, \{Q\in \PP(\C^{n\times n}): \innerp{A_j,Q}=0,\,\, j=1,\ldots,m \}.
\]
Note that $Z_{\mathcal A}$ is a linear subspace in $\PP(\C^{n\times n})$ with ${\rm dim}(Z_{\mathcal A}) \geq n^2-1-m$.
The projective variety $\MLP_{2r}^\C\subset \PP(\C^{n\times n})$ has dimension $4nr-4r^2-1$.
If $m\leq 4nr-4r^2-1$, then
\[
{\rm dim}(Z_{\mathcal A})+{\rm dim}(\MLP_{2r}^\C) \,\, \geq \,\, n^2-1,
\]
which implies that (see \cite[Prop.11.4]{alge})
\[
Z_{\mathcal A} \cap \MLP_{2r}^\C \neq \emptyset.
\]
Hence, if $m\leq 4nr-4r^2-1$  there exits a non-zero  matrix $Q_0\in Z_{\mathcal A} \cap \MLP_{2r}^\C $ satisfying
\[
\innerp{A_j,Q_0}=0,\quad\quad j=1,\ldots, m
\]
which implies ${\bf M}_{\mathcal A}$ is not injective on $\ML_r^\C$.
\end{proof}
\begin{remark}
We also can use the technology in the proof of Theorem \ref{th:comp} to study the {\em weak recovery }, which means to
recover {\em a fixed } $Q_0\in \ML_r^\C$ from ${\bf M}_{\mathcal A}(Q_0)$ (see also \cite{uniq}). As shown in \cite{uniq}, to ensure ${\bf M}_{\mathcal A}$ has the weak recovery property,
we only need show that $Q=Q_0$ if ${\bf M}_{\mathcal A}(Q-Q_0)=0$ and $Q\in \ML_r^{\C}$. Note that
\[
\{Q-Q_0:Q\in \ML_r^\C\}\subset \C^{n\times n}
\]
is an affine variety with dimension $2nr-r^2$. Then using a similar method with the proof of Theorem \ref{th:comp}, we can show that ${\bf M}_{\mathcal A}$ is injective on $\{Q-Q_0:Q\in \ML_r^\C\}\subset \C^{n\times n}$ if $m\geq 2nr-r^2+1$ and $A_1,\ldots,A_m$ are $m$ generic matrices.
\end{remark}

\begin{remark}
It will be very interesting to construct $m=4nr-4r^2$ deterministic matrices $A_1,\ldots,A_m$ so that ${\bf M}_{\mathcal A}$ is injective on $\ML_r^\C$. In the context of phase retrieval, such constructions are presented in \cite{bodmann} and \cite{Philipp}. In fact, \cite{bodmann} and \cite{Philipp} present $4n-4$ deterministic Hermite matrices and prove they define an injective measurement on $\SS_1^\C$. It will be interesting to extend the results and methods of \cite{bodmann} and \cite{Philipp} to low-rank matrices  recovery.
\end{remark}
\section{The recovery of real low rank matrices}
In this section, we consider the case where $\H=\R$.
Then we have
\begin{theorem}\label{eq:realm}
Suppose that $r\leq n/2$.
Consider $m$ matrices ${\mathcal A}=\{A_1,\ldots,A_m\}\subset \R^{n\times n}$ and the mapping ${\bf M}_{\mathcal A}:\R^{n\times n}\rightarrow \R^m$.  The following holds

(a) If $m\geq 4nr-4r^2$ then ${\bf M}_{\mathcal A}$ is injective on $\ML_r^\R$ for generic matrices $A_1,\ldots,A_m$.

(b) Suppose that $n=2^k+r, k\in\Z_+,$ or $n=2r+1$. If $m<4nr-4r^2$, then ${\bf M}_{\mathcal A}$ is  not injective on $\ML_r^\R$.

\end{theorem}
\begin{proof}
The proof of Part (a) is similar with the proof of (a) in Theorem \ref{th:comp} and hence we omit it here. We next turn to (b).
Following the notation from the proof of Theorem \ref{th:comp}, we set
\[
Z_{\mathcal A} \,\,:=\,\, \{Q\in \PP(\C^{n\times n}): \innerp{A_j,Q}=0,\,\, j=1,\ldots,m \}.
\]
Note that $Z_{\mathcal A}$ is a linear space and ${\rm dim}(Z_{\mathcal A}) \geq n^2-1-m$.
To state conveniently, set
\[
\MLP_{2r}^\C:=\{Q\in \PP(\C^{n\times n}): {\rm rank}(Q)\leq 2r\}.
\]
The $\MLP_{2r}^\C$ is a projective variety in $\C^{n^2}$ and ${\rm dim}(\MLP_{2r}^\C) \geq 4nr-4r^2-1$.
Note that when $m\leq 4nr-4r^2-1$,
\[
{\rm dim}Z_{\mathcal A} + {\rm dim}(\MLP_{2r}^\C) \geq n^2-1,
\]
which implies that $\MLP_{2r}^\C\cap Z_{\mathcal A}\neq \emptyset$ \cite[Prop.11.4]{alge}. According to Lemma \ref{le:odd}, the variety $\MLP_{2r}^\C$ has odd degree provided $n=2^k+r$. Note that $Z_{\mathcal A}$ is a linear space hence the intersection between
$\MLP_{2r}^\C$ and $Z_{\mathcal A}$ also has odd degree, which implies that the intersection
$\MLP_{2r}^\C\cap Z_{\mathcal A}$ has a real point since any projective variety with odd degree defined over $\R$ has real point. Thus
 there exists a nonzero real matrix $Q_0\in \MLP_{2r}^\C\cap Z_{\mathcal A}$, which implies that ${\bf M}_{\mathcal A}$ is  not injective on $\ML_r^\R$.

\end{proof}

According to \cite[Ex. 19.10]{alge}, the degree of the projective variety of $\MLP_{2r}^\C$ is
\[
d_{n,2r}:=\prod_{i=0}^{n-2r-1}\frac{(n+i)!\cdot i!}{(2r+i)!\cdot (n-2r+i)!}.
\]
Then
\begin{lem}\label{le:odd}
For $n=2^k+r, k\in \Z_+$ or $n=2r+1$, $d_{n,2r}$ is an odd integer.
\end{lem}
\begin{proof}
We first consider the case where $n=2r+1$.
A simple calculation shows that
\[
d_{2r+1,2r}=\frac{(2r+1)!}{(2r)!}=2r+1,
\]
which implies that $d_{n,2r}$ is odd provided $n=2r+1$.

We next assume that $n=2^k+r, k\in \Z_+$.
Note that
\begin{align*}
d_{n,2r}&=\prod_{i=0}^{n-2r-1}\frac{(n+i)!\cdot i!}{(2r+i)!\cdot (n-2r+i)!}=\prod_{i=0}^{n-2r-1}\frac{(n+i)\cdots  (n+i-(2r-1))}{(i+1)\cdots(i+2r)}\\
&= \prod_{i=0}^{2^k-r-1}\frac{(2^k+i+r)\cdots  (2^k+i-r+1)}{(i+1)\cdots(i+2r)}.
\end{align*}
Here, in the last equality, we use the assumption of $n=2^k+r$.
To state conveniently, we use $p_2(m)$ to denote the highest power of $2$ dividing $m\in \Z$. Then
\begin{align}
p_2(d_{n,2r})&=\sum_{i=0}^{2^k-r-1} \left(\sum_{j=-(r-1)}^r p_2(2^k+i+j)-\sum_{j=1}^{2r}p_2(i+j)\right)\nonumber\\
&=\sum_{i=0}^{2^k-r-1} \left(  \sum_{j=1}^r p_2(2^k+i+1-j)-\sum_{j=1}^{r}p_2(i+j+r) \right).\label{eq:p2}
\end{align}
Here, in the last equality, we use the fact of $p_2(i+j)=p_2(2^k+i+j)$ provided $i+j\leq 2^k-1$. We first consider the first term in (\ref{eq:p2}), i.e.,
\begin{align}
\sum_{i=0}^{2^k-r-1}  \sum_{j=1}^r p_2(2^k+i+1-j)&=\sum_{i=0}^{2r-1}  \sum_{j=1}^r p_2(2^k+i+1-j)+\sum_{i=0}^{2^k-3r-1}  \sum_{j=1}^r p_2(2^k+i+2r+1-j)\nonumber\\
&=\sum_{i=0}^{2r-1}  \sum_{j=1}^r p_2(2^k+i+1-j)+\sum_{i=0}^{2^k-3r-1}  \sum_{j=1}^r p_2(i+2r+1-j)\nonumber\\
&= \sum_{i=0}^{2r-1}  \sum_{j=1}^r p_2(2^k+i+1-j)+\sum_{i=0}^{2^k-3r-1}  \sum_{j=1}^r p_2(i+j+r).\label{eq:p3}
\end{align}
Then, combining (\ref{eq:p2}) and (\ref{eq:p3}), we obtain that
\begin{align*}
p_2(d_{n,2r})&= \sum_{i=0}^{2r-1}  \sum_{j=1}^r p_2(2^k+i+1-j)-\sum_{i=2^k-3r}^{2^k-r-1}  \sum_{j=1}^r p_2(i+j+r)\\
&= \sum_{i=0}^{2r-1}  \sum_{j=1}^r p_2(2^k+i+1-j)-\sum_{i=0}^{2r-1}\sum_{j=1}^r p_2(2^k+i+j-2r)\\
&= \sum_{i=0}^{2r-1}  \sum_{j=1}^r p_2(2^k+i+j-r)-\sum_{i=0}^{2r-1}\sum_{j=1}^r p_2(2^k+i+j-2r)\\
&= \sum_{i=r}^{2r-1}  \sum_{j=1}^r p_2(2^k+i+j-r)-\sum_{i=0}^{r-1}\sum_{j=1}^r p_2(2^k+i+j-2r)\\
&= \sum_{i=0}^{r-1}  \sum_{j=1}^r p_2(2^k+i+j)-\sum_{i=0}^{r-1}\sum_{j=1}^r p_2(2^k+i+j-2r)\\
&=\sum_{s=1}^{2r-1}(b_sp_2(2^k+s)-b_{2r-s}p_2(2^k-s))=0.
\end{align*}
Here, in the last equality, $b_s:=\#\{(i,j):i+j=s, 0\leq i \leq r-1, 1\leq j\leq r\}$ and
 we use the fact of $b_s=b_{2r-s}$ and $p_2(2^k+s)=p_2(2^k-s)$.
\end{proof}

Theorem \ref{eq:realm} shows the bound $4nr-4r^2$ is tight provided $n$ is in the form of $2^k+r$ or $2r+1$.
One may be interested in whether the bound $4nr-4r^2$ is tight in general.
The next counterexample shows that for the case where $(n,r)=(4,1)$ there exist $11=4n-5$ matrices ${\mathcal A}=\{A_1,\ldots,A_{11}\}\subset \R^{4\times 4}$ so that ${\bf M}_{\mathcal A}$ is  injective on $\ML_1^\R\subset \R^{4\times 4}$. And hence the bound is not tight provided $n=4, r=1$. We list the $11$ matrices as follows which are
obtained by computer random search:
{\tiny
\begin{equation}\label{eq:11matr}
  \begin{aligned}
A_1&= \begin{pmatrix} -4 & 1 & 3& 4 \\
 -4 & 4 & 4 & 3\\
 4 & -3 & 0 & -3 \\
  0 & -4 & 2 & 1 \\
 \end{pmatrix}
  A_2=\begin{pmatrix}
   0 & 3 & -1& -1 \\
 0 & -2 & -1 & 2\\
 0 & 3 & -2 & 3 \\
  1 & -1 & -3 & 2 \\
 \end{pmatrix}
 A_3= \begin{pmatrix}
   -1 & -4 & -1& -1 \\
 4 & 0 & -1 & 1\\
 -2 & 0 & 0 & 2 \\
  0 & -1 & 2 & 2 \\
 \end{pmatrix}
  A_4= \begin{pmatrix}
   -2 & -2 & 4& 1 \\
 -2 & 0 & 2 & 3\\
 1 & -2 & -4 & 3 \\
  -3 & 3 & 4 & -2 \\
 \end{pmatrix}\\
  A_5&=  \begin{pmatrix}
  4  & 2 & -4& -4 \\
 -4 & -3 & 0 & 0\\
 1 & -4 & 4 & -2 \\
  3 & 0 & 2 & 0 \\
 \end{pmatrix}
A_6= \begin{pmatrix}
 2  & 2 & 3& 4 \\
 2 & -4 & 3 & 1\\
 0 & -2 & 1 & -2 \\
  -1 & 0 & -1 & -4 \\
 \end{pmatrix}
 A_7=\begin{pmatrix}
 2  & 1 & 4& 0 \\
-1 & -3 & 0 & -1\\
 4 & -1 & -4 & 3 \\
 0 & 3 & 0 & 4 \\
 \end{pmatrix}
 A_8=\begin{pmatrix}
 0  & 3 & -1& 2 \\
4 & 2 & 1 & 1\\
 -2 & -1 & 3 & 4 \\
 3 & 0 & 3 & 3 \\
 \end{pmatrix}\\
  A_9&=\begin{pmatrix}
2  & -1 & 4& -4 \\
-2 & 2 & 3 & -1\\
 -1 & 1 & 4 & -1 \\
 -3 & -4 & 4 & 3 \\
 \end{pmatrix}
  A_{10}= \begin{pmatrix}
-4  & 2 & 0& -1 \\
4 & 1 & 0 & 4\\
 -1 & -3 & 4 & 1 \\
 -3 & 2 & 4 & -4 \\
 \end{pmatrix}
   A_{11}= \begin{pmatrix}
1  & 1 & -2& 0 \\
3 & 0 & -2 & -4\\
 2 & -4 & -2 & 4 \\
 4 & 3 & 2 & -2 \\
 \end{pmatrix}.
 \end{aligned}
\end{equation}}
 We show the map ${\bf M}_{\mathcal A}$ associated with the $11$ matrices is injective on $\ML_1^\R$:
\begin{theorem}\label{th:cont}
Set ${\mathcal A}=\{A_1,\ldots,A_{11}\}$ where $A_{j}, j=1,\ldots,11,$ are defined in (\ref{eq:11matr}). Then the map ${\bf M}_{\mathcal A}$
is injective on $ \ML_1^\R\subset \R^{4\times 4}$.
\end{theorem}
\begin{proof}
To this end, we only need prove the set
\begin{align}\label{eq:Q}
\{Q\in \R^{4\times 4}: \innerp{A_j,Q}=0, j=1,\ldots, 11, {\rm rank}(Q)\leq 2\}
\end{align}
has only zero matrix.

We build the proof following the ideas of  Vinzant \cite[Theorem 1]{small}.
In fact, we use Vinzant's test, which is stated in Algorithm \ref{alg:theorem1}, to verify the map ${\bf M}_{\mathcal A}$
is injective on $ \ML_1^\R\subset \R^{4\times 4}$. We next explain the reason why Algorithm \ref{alg:theorem1} works.
 Any $4\times 4$ real matrix can be written as
\[
 Q= \begin{pmatrix}
   x_{11} & x_{12} & x_{13}& x_{14} \\
 x_{21} & x_{22} & x_{23} & x_{24}\\
 x_{31} & x_{32} & x_{33} & x_{34} \\
  x_{41} & x_{42} & x_{43} & x_{44} \\
 \end{pmatrix},
\]
where $x_{jk}, 1\leq j\leq 4, 1\leq k\leq 4,$ are $16$  variables.  Set
\[
\ell_j:=\innerp{A_j, Q},\,\, j=1,\ldots,11
\]
and we use $m_{jk}$ to denote the determinant of the sub-matrix  formed by deleting the $j$th row and $k$th column from the matrix $Q$. Note that both $\ell_j$ and $m_{jk}$ are  polynomials about $x_{11},\ldots, x_{44}$.
We recall the fact ${\rm rank}(Q)\leq 2$ is equivalent to the vanish of $m_{jk}, j=1,\ldots,4, k=1,\ldots,4$. Hence, he map ${\bf M}_{\mathcal A}$
is injective if and only if the polynomial system
\begin{equation}\label{eq:poly}
m_{11}=m_{12}=\cdots=m_{44}=\ell_1=\cdots=\ell_{11}=0
\end{equation}
has no nonzero real solution $(x_{11},\ldots,x_{44})\in \R^{16}$. A simple observation is that if ${\bf x}^0:=(x_{11}^0, x_{12}^0,\ldots,x_{44}^0)\in \R^{16}$ is a root of (\ref{eq:poly}) then $f({\bf x}^0)=0$ for any $f$ in the ideal generated by the set of  polynomials $\{ m_{11},\ldots,m_{44},\ell_{1},\ldots,\ell_{11}\}$. To state conveniently, we use the notation $\left< m_{11},\ldots,m_{44},\ell_{1},\ldots,\ell_{11} \right>$ to denote the ideal. We use the computer algebra software {\tt maple} to compute a Gr\"{o}bner basis of the ideal   and elimination.
The result is a polynomial $f_0\in \Q[x_{43},x_{44}]$ (see Appendix A), which is a homogeneous polynomial of degree $20$. Then $f_0(x_{43}^0,x_{44}^0)=0$ if ${\bf x}^0=(x_{11}^0, x_{12}^0,\ldots,x_{44}^0)\in \R^{16}$ is a root of (\ref{eq:poly}) since $f_0\in \left< m_{11},\ldots,m_{44},\ell_{1},\ldots,\ell_{11} \right>$. We claim that $f_0$ has only real root $(0,0)$. Indeed, if $f_0$ has a nonzero real solution $(x_{43}^0,x_{44}^0)$ then $x_{44}^0\neq 0$ (otherwise, $x_{43}^0=x_{44}^0=0$ since $f_0$ is a homogeneous polynomial). Note that if $f_0(x_{43}^0,x_{44}^0)=0$ then $f_0(\lambda x_{43}^0,\lambda x_{44}^0)=0$ for any $\lambda \in \C$. Without loss of generality, we suppose that $x_{44}^0=1$. Then using Sturm  sequences, we can verify the univariate polynomial $f_0(x_{43},1)$ has no real solutions, which implies that the real roof of $f_0(x_{43},x_{44})$ is only $(0,0)$.

We claim that there is nonzero root to (\ref{eq:poly}) with $x_{43}=0, x_{44}=0$. We can verify the claim still by Gr\"{o}bner basis. For any $\lambda\in \C$, $(\lambda x_{11}^0, \lambda x_{12}^0,\ldots,\lambda x_{44}^0)$ is  a solution to (\ref{eq:poly}) if $(x_{11}^0, x_{12}^0,\ldots, x_{44}^0)$ is a root to   (\ref{eq:poly}). Hence if (\ref{eq:poly}) has a nonzero solution, then (\ref{eq:poly}) will have a solution with some coordinate  equal to 1. We first consider the case where $x_{11}=1$. By computing a Gr\"{o}bner basis of the ideal $\left<x_{11}-1, x_{43},x_{44}, m_{11},\ldots,m_{44},\ell_1,\ldots,\ell_{11}\right>$, we can verify
\[
1\in \left<x_{11}-1, x_{43},x_{44}, m_{11},\ldots,m_{44},\ell_1,\ldots,\ell_{11}\right>,
\]
which implies that (\ref{eq:poly}) has no solution with $x_{11}=1, x_{43}=0$ and $x_{44}=0$. We can verify the case $x_{jk}=1$ with $(j,k)\in [1,4]^2\setminus \{(4,3), (4,4)\}$ one by one.

 We post the code for these computation in {\tt Maple }  at http://lsec.cc.ac.cn/$\sim$xuzq/LowRank.htm.
\end{proof}

\begin{algorithm}[htb]
\caption{
Vinzant's test for injective of the map  ${\bf M}_{\mathcal A}$ }
\label{alg:theorem1}
\begin{algorithmic}[1]
\REQUIRE m=11, the matrices  $A_1,\ldots, A_{m}$ which are given in (\ref{eq:11matr}),
   \[
 Q= \begin{pmatrix}
     x_{11} & x_{12} & x_{13}& x_{14} \\
 x_{21} & x_{22} & x_{23} & x_{24}\\
 x_{31} & x_{32} & x_{33} & x_{34} \\
  x_{41} & x_{42} & x_{43} & x_{44} \\
 \end{pmatrix}.
\]\par
\STATE Set $Q_{jk},1\leq j,k\leq 4$ is the sub-matrix of $Q$ formed by deleting  $j$th column and $k$th row from $Q$.
\STATE Set $\ell_j=\innerp{A_j,Q}, j=1,\ldots,m$ and $m_{jk}=\det(Q_{jk}), 1\leq j,k\leq 4$.
\STATE Computer Gr\"{o}bner basis of the ideal $\left<\ell_1,\ldots,\ell_{m}, m_{11},m_{12},\ldots,m_{44}\right>$
and obtain that $f_0\in \Q[x_{4,3},x_{4,4}]$.
\STATE Use Sturm Sequence to compute the number of real roots of $f_0(x_{4,3},1)$.
\IF{the number of real roots of $f_0(x_{4,3},1)$ is $0$ }
\FORALL{ $j,k\in [1,4]\times [1,4]$  }
\STATE Check whether  $1\in \left<x_{j,k}-1, x_{4,3},x_{4,4},\ell_1,\ldots,\ell_{m},m_{11},\ldots,m_{44}\right>$ by computing Gr\"{o}bner basis
\IF {$1\in \left<x_{j,k}-1, x_{4,3},x_{4,4},\ell_1,\ldots,\ell_{m},m_{11},\ldots,m_{44}\right>$}
\STATE $r_{j,k}=1$
\ELSE
\STATE $r_{j,k}=0$, ``FAIL"
\ENDIF
\IF {$r_{j,k}=1$ for all $j,k\in [1,4]\times [1,4]$ }
\STATE ``INJECTIVE"
\ENDIF
\ENDFOR
\ELSE
\STATE ``FAIL"
\ENDIF
\end{algorithmic}
\end{algorithm}

\subsection{Symmetric matrix}

We next consider the symmetric matrix which will be helpful for the investigation of phase retrieval by projection. Recall that
\[
\SS_{r}^\R:=\{X\in \R^{n\times n}: {\rm rank}(X)\leq r, X^\top=X\}.
\]

\begin{theorem}\label{th:symm}
Suppose that $r\leq n/2$.
Consider $m$ matrices $A_1,\ldots,A_m\in  \R^{n\times n}$ and the mapping ${\bf M}_{\mathcal A}:\R^{n\times n}\rightarrow \R^m$
where ${\mathcal A}=\{A_1,\ldots,A_m\}$.  The following holds

(a) If $m\geq 2nr+r-2r^2$ then ${\bf M}_{\mathcal A}$ is injective on $\SS_r^\R$ for generic matrices $A_1,\ldots,A_m$.

(b) If $n=2^k+r$ and $m<2nr+r-2r^2$, then ${\bf M}_{\mathcal A}$ is  not injective on $\SS_r^\R$.
\end{theorem}
\begin{proof}
A simple observation is that the map ${\bf M}_{\mathcal A}$ is injective on $\SS_r^\R$ if and only if there is a nonzero $Q\in \SS_{2r}^\R$ for which
 ${\bf M}_{\mathcal A}(Q)=0$.
Thus, we only need show that $Q=0$ provided $Q\in \SS_{2r}^\R$ and  ${\bf M}_{\mathcal A}(Q)=0$.
 Recall that $\SS_r^\R$ is an affine algebraic variety with dimension ${n+1\choose 2}-{n-r+1\choose 2}$, which implies that
${\rm dim}(\SS_{2r}^\R)=2nr+r-2r^2$.
The proof of Part (a) is similar with the proof of (a) in Theorem \ref{th:comp} and we omit it here.
 We next turn to (b).
We set
\begin{align*}
Z_{\mathcal A} \,\,&:=\,\, \{Q\in \PP(\C^{n\times n}): \innerp{A_j,Q}=0, Q^\top=Q,\,\, j=1,\ldots,m \},\\
\mathcal{PS}_{2r}^\C&:=\{X\in \PP(\C^{n\times n}): {\rm rank}(X)\leq 2r,X^\top=X\}.
\end{align*}
Then ${\rm dim}(Z_{\mathcal A}) \geq \frac{n(n+1)}{2}-1-m$
and ${\rm dim}(\mathcal{PS}_{2r}^\C) \geq 2nr+r-2r^2-1$.
Note that when $m\leq 2nr+r-2r^2-1$,
\[
{\rm dim}Z_{\mathcal A} +{\rm dim}(\mathcal{PS}_{2r}^\C)  \geq \frac{n(n+1)}{2}-1,
\]
which implies that $\mathcal{PS}_{2r}^\C \cap Z_{\mathcal A}\neq \emptyset$ \cite[Prop.11.4]{alge}.
 Note that $ \mathcal{PS}_{2r}^\C =\MLP_{2r}^\C \bigcap_{1\leq j,k\leq m}{\mathcal H}_{jk}$,
 where
 \[
 {\mathcal H}_{jk}:=\{X\in \PP( \C^{n\times n}): x_{jk}=x_{kj}\}.
 \]
 According to Lemma \ref{le:odd}, the variety $\MLP_{2r}^\C$ has odd degree provided $n=2^k+r$, which implies that
the degree of $ \mathcal{PS}_{2r}^\C$ is odd if $n=2^k+r$ since  $\mathcal{PS}_{2r}^\C$ is the intersection of $\MLP_{2r}^\C$ and some linear spaces. Note that $Z_{\mathcal A}$ is a linear space hence the intersection between
$\mathcal{PS}_{2r}^\C$ and $Z_{\mathcal A}$ also has odd degree, which implies that the intersection
$\mathcal{PS}_{2r}^\C\cap Z_{\mathcal A}$ has a real point since any projective variety with odd degree defined over $\R$ has real point.

\end{proof}

\begin{remark}
When $r=1$, the bound $2nr+r-2r^2$ is reduced to $2n-1$. A natural question is whether the bound $2n-1$ is tight for the recovery
of the symmetric rank-1  matrix.
For the case $n=4$, one can construct $6=2n-2$ matrices which are injective on $\mathcal{PS}_{1}^\R$ (see Theorem \ref{th:42}), which implies that
the bound $2nr+r-2r^2$ is not tight for $n=4,r=1$.
\end{remark}

\begin{remark}
If we require  $A_j$ is in the form of $a_ja_j^\top$ with $a_j\in \R^n$, then the bound $2n-1$ is tight. In fact, $\innerp{A_j, Q}=\abs{\innerp{a_j,x}}^2$ provided
$Q=xx^\top\in \R^{n\times n}$.  According to the result from phase retrieval \cite{BCE}, ${\bf M}_{\mathcal A}$ is injective on $\SS_1^\R$ if and only if $\{a_1,\ldots,a_m\}\subset \R^n $ satisfies the finite complement property, i.e., for every subset $I\subset\{1,\ldots,m\}$ either $\{a_j\}_{j\in I}$ or $\{a_j\}_{j\in I^c}$ spans $\R^n$, which implies the bound $2n-1$ is tight provided the measurement matrices $A_j$  is in the form of $a_ja_j^\top$.
\end{remark}

\section{phase retrieval by projections}

Recall that we use  $P_j:\R^n\rightarrow W_j$ to denote an  orthogonal projection where $W_j\subset \R^n$ is a subspace.
As introduced in  Section 1.3, we say that $\{W_j\}_{j=1}^m$ yields phase retrieval if for all $x,y\in \R^n$ satisfying $\|P_jx\|=\|P_jy\|$ for all $j=1,\ldots,m$ then $x=\pm y$.
The following theorem shows that $2n-1$ projections are enough for phase retrieval by projection.
\begin{theorem}\cite{phaseproj}
Phase retrieval can be done in $\R^n$ with $2n-1$ subspaces each of any dimension less than $n-1$.
\end{theorem}
 The problem is also raised in \cite{phaseproj1} which states {\em can phase
retrieval be done in $\R^n$ with fewer than $2n-1$ projections?}  Based on Theorem \ref{th:symm}, we show that the bound $2n-1$ is tight provided $n=2^k+1, k\in \Z_{\geq 1}$. Particularly, we have
\begin{coro}
Suppose that $n$ is in the form of $2^k+1$.
Given any subspaces $\{W_j\}_{j=1}^m$ in $\R^n$ with $m<2n-1$, there exist $x,y\in \R^n$ with $x\neq \pm y$ so that
$\|P_jx\|=\|P_jy\|, j=1,\ldots,m$.
\end{coro}
\begin{proof}
Suppose that $\{u_{j,1},\ldots,u_{j,d_j}\}$ is an orthonormal basis of $W_j$. Then for any $x\in \R^n$
\[
\|P_jx\|^2=\innerp{A_j,xx^\top},\quad j=1,\ldots,m
\]
where
\begin{equation}\label{eq:Aj}
A_j:=u_{j,1}u_{j,1}^\top+\cdots+u_{j,d_j}u_{j,d_j}^\top.
\end{equation}
 Then $\{W_j\}_{j=1}^m$ allow phase retrieval if and only if
the map ${\bf M}_{\mathcal A}$ is injective on $\SS_1^\R$ where ${\mathcal A}=\{A_1,\ldots,A_m\}$ and $A_j, j=1,\ldots,m$ are defined in (\ref{eq:Aj}).  The part (b) in Theorem \ref{th:symm}
shows that ${\bf M}_{\mathcal A}$ is not injective if $m<2n-1$ and $n=2^k+1$ which implies the corollary.
\end{proof}

Naturally, one may be interested in whether the bound $2n-1$ is tight when $n\neq 2^k+1$.
We give a negative answer by  presenting a counterexample for the case where $n=4$. In fact, we present $7$ subspaces in $\R^4$, which
is obtained by the computer search.
Set
\begin{align}\label{eq:UV}
  U:=[u_1,u_2,\ldots,u_6]= \begin{pmatrix}
   1 & 0 & 0& 0 & 1 & 1 \\
 0 & 1 & 0 & 0 &0 & 1\\
 0 & 0 & 1 & 0 &-4 & 2\\
  0 & 0 & 0 & 1 &-3& 3\\
 \end{pmatrix}
V:=[v_1,v_2,\ldots,v_6]=
 \begin{pmatrix}
   0 & 5 & -1& -1 & -17 & -5 \\
 -5 & 0 & 0 & 5 &4 & 4\\
 2 & -2 & 0 & -3 &-2 & 2\\
  -2 & 1 & 0 & 0 &-3& -1\\
 \end{pmatrix}
\end{align}
and
\begin{equation}\label{eq:W}
W_j:={\rm span}\{u_j,v_j\}\subset \R^4,\quad j=1,\ldots,6.
\end{equation}

\begin{theorem}\label{th:42}
Suppose that $W_1,\ldots,W_6$ are defined in (\ref{eq:W}).  Then the phase retrieval by projection can be done in $\R^4$ with
the 6 subspaces $W_1,\ldots,W_6$.
\end{theorem}
\begin{proof}
To this end, we only need show that the set
\[
\{Q\in \R^{4\times 4}: {\rm rank}(Q)\leq 2,\, Q^\top=Q,\, \innerp{A_j,Q}=0,\, j=1,\ldots,6\},
\]
has only zero matrix,
where
\begin{equation}\label{eq:Ajuv}
A_j=\frac{1}{\|u_j\|^2}u_ju_j^\top+\frac{1}{\|v_j\|^2}v_jv_j^\top, j=1,\ldots,6.
\end{equation}
Any symmetric $4\times 4$ matrix can be written as
\begin{equation}\label{eq:Q}
 Q= \begin{pmatrix}
   x_{11} & x_{12} & x_{13}& x_{14} \\
 x_{12} & x_{22} & x_{23} & x_{24}\\
 x_{13} & x_{23} & x_{33} & x_{34} \\
  x_{14} & x_{24} & x_{34} & x_{44} \\
 \end{pmatrix}.
\end{equation}
Then we can verify ${\bf M}_{\mathcal A}$ is injective on $\SS_1^\R$ with ${\mathcal A}=\{A_1,\ldots,A_6\}$ by using a similar method with
Algorithm \ref{alg:theorem1}. In fact, we
 verify ${\bf M}_{\mathcal A}$ is injective  by Algorithm 1
with inputting  $m=6$, $A_1,\ldots, A_m$, and $Q$, which are given in (\ref{eq:Ajuv}) and  (\ref{eq:Q}), respectively. In Line 3 of Algorithm 1,
we obtain $f_0\in \Q[x_{34},x_{44}]$ by computing the Gr\"{o}bner basis, which is shown as follows:

{\tiny
\begin{align*}
f_0(x_{3,4},x_{44})=&519966562263643554945384191703616165395119112637573956248783182163623419737152512\cdot x_{3,4}^{10}\\
&+328579249789044588040378180884308246920612283125532941525689039614516382331655616\cdot x_{3, 4}^9x_{4,4}\\
&-1488937659336445964244382640314269130820570919591539086664937634093065289567722168\cdot x_{3, 4}^8x_{4,4}^2\\
&-1233940048680917718690405606336150050552321029972511785557378326694535160186545716\cdot x_{3, 4}^7x_{4,4}^3\\
&+560830252887171704842938122056614428129006430576409737097486332815580249915632343\cdot x_{4, 4}^4x_{3, 4}^6\\
&+862820282355455834964156023668088334391251839891193308388567471700918807326409512\cdot x_{3, 4}^5x_{4,4}^5\\
&+775754320104988082038883422728415511210018447819392679061447643613437650180479290\cdot x_{3, 4}^4x_{4, 4}^6\\
&+470277354383117587803315722906687832056957450220672159363137099743786021781725000\cdot x_{3, 4}^3x_{4, 4}^7\\
&+200124465435786576051259835347973959273830283051981724273474893893577713039639000\cdot x_{3, 4}^2x_{4, 4}^8\\
&+61607571777035859344852037093432432693206774944581030976959966534882027371910000\cdot x_{3, 4}x_{4, 4}^9\\
&+8636626929016108140668606241999544716256996255890976455134579423688060066275000\cdot x_{4, 4}^{10}.
\end{align*}
}
We post the code for these computation in {\tt Maple }  at http://lsec.cc.ac.cn/$\sim$xuzq/LowRank.htm.
\end{proof}

\begin{appendix}

\begin{landscape}
\section{The $f_0(x_{43},x_{44})$ which is used in the proof of Theorem \ref{th:symm}}
{\tiny
\begin{align*}
& 2190263004585315683793318506387418373659624458975581329201799757712208390076105396569258443797012006452503150042124624662199923225090229659511210975239970888153158\cdot x_{43}^{20}\\
& -5179707489822384879036242558259786010072295064540826084187984072746733580158433711332940236665034062027211264655232956957554594817
 844803752968192178160111903891202\cdot x_{43}^{19}\cdot x_{44}\\
& +1040712214012933805249563659200424890987717333211994380167449996686920813620554097934268846901580572027138626281857992010765011147935
 72262681595311459946577333645090\cdot x_{44}^2\cdot x_{43}^{18}\\
& -216072112689014792523323463032129738500684274319297996378543592789980269601501663337459487098688654806444070799488716026096089429093
895413888228401366614027073502628\cdot x_{43}^{17}\cdot x_{44}^3\\
& +150656214432848705569353420318391110491666871524757876108892401746165136597521350258275954384514944040061546271208263265306517607658572682
2027093617595286061839632576\cdot x_{44}^4\cdot x_{43}^{16}\\
&+299875005069982324464737232616434147573270672330608059881367750159691634086066489001056800936147027757951078101083160502892144438351541377
231003456448761936746488414\cdot x_{43}^{15}\cdot x_{44}^5\\
&-1636431767115105327642256045480119509672356313030010174315380376677183306233846137064478853453247718135097465652724300959059943528386729354809058562415946730020759127\cdot x_{43}^{14}\cdot x_{44}^6\\
&+70489634113319253768763757912759898026710660340641300929928752453490269002171013933930144336191577512046605593483286614555698346489612902425503079679048677581991711201\cdot x_{43}^{13}\cdot x_{44}^7\\
&-2837512873942817681084742778680818174664543414063278928066522307542932531478416818140806071570393059086601600717697771957536549044543060978
00824447679379474504152094942\cdot x_{43}^{12}\cdot x_{44}^8\\
&+995932712734255470529235127196212532582344056670644822018459702834698228411481231741127911879071226987664884637125259887532399285303825908296757977886041651714696917778\cdot x_{43}^{11}\cdot x_{44}^9\\
&-258445808607953313164679652132032036429715368478853009379862012348067562469944501088340890453492422913568007790770281431624985291845338328430
8052222356363717915656742284\cdot x_{43}^{10}\cdot x_{44}^{10}\\
&+5722961393302102292758678488211850165342295398495413012439722250125847472613133015046362267108099705458400660879775196934139225875722766167526818738883868235472207295318\cdot x_{43}^9\cdot x_{44}^{11}\\
&-630903666163102044269305678828872932653200564493312344178612822228642037237561527877847065585168126464704141871730257087486167130917377572099
708190266736488397171274544\cdot x_{43}^8\cdot x_{44}^{12}\\
&-66044655033880585342698995668204656740186018128340889546801734092786485152939244784658016013861989847047384410257733754028999986653323322558108522264030232060720866849482\cdot x_{43}^7\cdot x_{44}^{13}\\
&+331710411979734991930045202128494805001764789107926958310437527724295580780866217730043775840527356748728419672335038312502974307254601559883712666369538857008427480797489\cdot x_{43}^6\cdot x_{44}^{14}\\
&-965845696733814029786999463978506859651818142902503987920109518394356204424235798108650852569849782624138330160838049982847845141054922144933695186034515159397327737110595\cdot x_{43}^5\cdot x_{44}^{15}\\
&+1910472054626214661171014521853950396829436381820503787956311058855103556898420463438512319957613171396422456016287996033426353168337725276189486781016250884277767002962296\cdot x_{43}^4*x_{44}^{16}\\
&-2572438761515126654331167975391717741834396831399029211381124450074460431446556299914852071543732665803310127498229219581211443279655430825498112043444924243569330420969108\cdot x_{43}^3\cdot x_{44}^{17}\\
&+2317118396824473116957993520633692645043414501959594440441532677816205844606801169352009501563550170576981690146529499525917896571663586948450976618053393044056987915792096\cdot x_{43}^2\cdot x_{44}^{18}\\
&-1266178703505322402605611812456953754907639982077075593911044325271865979988606348124211651128143828443791667351751806483036900444172736400483269837445165843264333005502240\cdot x_{44}^{19}\cdot x_{43}\\
&+341017021380913068697699063908437633507458704746930441702955178001151569680017455954493009794352676756966784216718914427692248083290003365728420621063417637683002937388800\cdot x_{44}^{20}.
\end{align*}}
\end{landscape}
\end{appendix}


\begin{thebibliography}{99}


\bibitem{phaseproj1}Saeid Bahmanpour, Jameson Cahill, Peter G. Casazza, John Jasper, Lindsey M. Woodland, Phase retrieval and norm retrieval,
 arXiv:1409.8266.

\bibitem{BCE} Radu Balan, Pete Casazza and Dan Edidin, On signal reconstruction without phase, Applied and Computational Harmonic Analysis
Volume 20, Issue 3, May 2006, Pages 345¨C356.



\bibitem{BCMN} Afonso S. Bandeira, Jameson Cahill, Dustin G. Mixon,  Aaron A. Nelson,
Saving phase: Injectivity and stability for phase retrieval,
Applied and Computational Harmonic Analysis
Volume 37, Issue 1, July 2014, Pages 106-125

\bibitem{bodmann}
B. G. Bodmann and N. Hammen,
{\em Stable phase retrieval with low-redundancy frame}, Advances in Computational Mathematics
 41(2015), 317-331.


\bibitem{phaseproj}Jameson Cahill, Peter G. Casazza, Jesse Peterson, Lindsey Woodland, Phase Retrieval By Projections, arXiv:1305.6226.


\bibitem{LMshen}J. F. Cai, E.J. Cand\`{e}s, Z.W. Shen, A singular value thresholding algorithm for matrix completion, SIAM J. Optim., 20(4)(2008), 1956-1982.

\bibitem{LMPlan}    E. J. Cand\`{e}s, Y. Plan, Tight Oracle Inequalities for Low-Rank Matrix Recovery From a Minimal Number of Noisy Random Measurements, IEEE Trans. Inform. Theory 57(4)(2010)2342-2359.

\bibitem{LMCR}E. J. Cand\`{e}s and B. Recht, Exact matrix completion via convex optimization, Foundations of Computational Mathematics 9(6), 717-772 (2009).

\bibitem{ROP} T. Tony Cai, Anru Zhang,  ROP: Matrix recovery via rank-one projections, arXiv:1310.5791.


\bibitem{CEHV} Aldo Conca, Dan Edidin,  Milena Hering and Cynthia Vinzant,
An algebraic characterization of injectivity in phase retrieval, Applied and Computational Harmonic Analysis, Volume 38, Issue 2, March 2015, Pages 346-356.


\bibitem{uniq} Y. C. Eldar, D. Needell, Y. Plan, Uniqueness conditions for low-rank matrix recovery, Appl. Comp. Harm. Anal, 33(2013)309-314.
\bibitem{alge} Joe Harris, {\em Algebraic geometry } Vol. 133 of GTM, Springer-Verlag, New York, 1992. A first course, corrected reprint of the 1992 original.
\bibitem{KRT} Richard Kueng, Holger Rauhut, Ulrich Terstiege,  Low rank matrix recovery from rank one measurements, arXiv:1410.6913.


\bibitem{His2}    S. Oymak,  and B. Hassibi, New Null Space Results and Recovery Thresholds for Matrix Rank Minimization, arXiv:1011.6326.
\bibitem{His3}P. Parrilo, A. Willsky, V. Chandrasekaran, and B. Recht, The Convex Geometry of Linear Inverse Problems, Foundations of Computational Mathematics December 2012, Volume 12, Issue 6, pp 805-849.
\bibitem{Philipp}
Friedrich Philipp, Phase Retrieval from 4N-4 Measurements, arXiv:1403.4769.

\bibitem{LMRe} B. Recht, M. Fazel, and P. Parrilo, Guaranteed minimum-rank solutions of linear matrix equations via
nuclear norm minimization, SIAM Rev. 52(3)(2010)471-501.



\bibitem{His1} V.Y.F. Tan,  L. Balzano and    S.C. Draper,  Rank Minimization Over Finite Fields: Fundamental Limits and Coding-Theoretic Interpretations,  IEEE Transactions on  Information Theory, 58(4)(2012)2018-2039.

\bibitem{small}Cynthia Vinzant, A small frame and a certificate of its injectivity, arXiv:1502.04656.







\end{thebibliography}
\end{document}